\documentclass{amsart}
\usepackage[T1]{fontenc}
\usepackage{ae}
\usepackage{aecompl}%
\usepackage{ucs}
\usepackage[utf8x]{inputenc}
\usepackage{color}
\usepackage[enableskew,vcentermath]{youngtab}
\usepackage{amsfonts}
\usepackage{hhline}
\usepackage{amsthm}
\usepackage{amsmath}
\usepackage{amssymb}
\usepackage[all]{pstricks}
\usepackage[justification=centering,labelfont=bf]{caption}
\usepackage{hyperref}

\begin{document}

\numberwithin{equation}{section}

\newtheorem{Le}{Lemma}[section]
\newtheorem{Ko}[Le]{Lemma}
\newtheorem{Sa}[Le]{Theorem}
\newtheorem{pro}[Le]{Proposition}

\newtheorem{theorem}[Le]{Theorem}
\newtheorem{lemma}[Le]{Lemma}
\newtheorem{Con}[Le]{Conjecture}
\newtheoremstyle{Bemerkung}
  {}{}{}{}{\bfseries}{.}{0.5em}{{\thmname{#1}\thmnumber{ #2}\thmnote{ (#3)}}}

\theoremstyle{Bemerkung}
\newtheorem{definition}[Le]{Definition}
\newtheorem{Def}[Le]{Definition}
\newtheorem{example}[Le]{Example}

\newtheorem{remark}[Le]{Remark}
\newtheorem{Bem}[Le]{Remark}
\newtheorem{Bsp}[Le]{Example}

\renewcommand{\l}{\lambda}
\renewcommand{\L}{\lambda}
\newcommand{\bl}{\bar\lambda}
\newcommand{\bn}{\bar\nu}
\newcommand{\bg}{\bar g}
\newcommand{\A}{\mathcal{A}}
\newcommand{\B}{\mathcal{B}}
\newcommand{\C}{\mathcal{C}}
\newcommand{\D}{\mathcal{D}}
\renewcommand{\S}{\mathcal{S}}
\renewcommand{\L}{\mathcal{L}}
\newcommand{\E}{\mathcal{E}}
\newcommand{\N}{\mathcal{N}}
\renewcommand{\a}{\alpha}
\renewcommand{\b}{\beta}
\renewcommand{\c}{\gamma}
\renewcommand{\d}{\delta}
\newcommand{\e}{\epsilon}
\newcommand{\h}{\hfil}
\newcommand{\X}{X}
\newcommand{\abs}[1]{\left| #1 \right|}
\newcommand{\lm}{\l/\mu}
\renewcommand{\ln}{\l/\nu}
\newcommand{\ab}{\a/\b}
\newcommand{\m}{\mu}
\newcommand{\n}{\nu}
\newcommand{\bp}{\bar p}

\title[Reduced Kronecker Products]{Reduced Kronecker products which are multiplicity free or contain only few components}
\author[C. Gutschwager]{Christian Gutschwager}
\address{Institut für Algebra, Zahlentheorie und Diskrete Mathematik, Leibniz Universität Hannover,  Welfengarten 1, D-30167 Hannover}
\email{gutschwager (at) math (dot) uni-hannover (dot) de}


\begin{abstract}
It is known that the Kronecker coefficient of three partitions is a bounded and weakly increasing sequence if one increases the first part of all three partitions. Furthermore, if the first parts of partitions $\l,\m$ are big enough then the coefficients of the Kronecker product $[\l][\m]=\sum_\n g(\l,\m,\n)[\n]$ do not depend on the first part but only on the other parts. The reduced Kronecker product $[\l]_\bullet\star[\m]_\bullet$ can be viewed (roughly) as the Kronecker product $[(n-\abs{\l},\l)][(n-\abs{\m},\m)]$ for $n$ big enough.

In this paper we classify the reduced Kronecker products which are multiplicity free and those which contain less than $10$ components. Furthermore, we give general lower bounds for the number of constituents and components of a given reduced Kronecker product. We also give a lower bound for the number of pairs of components whose corresponding partitions differ by one box.

Finally we argue that equality of two reduced Kronecker products is only possible in the trivial case that the factors of the product are the same.
\end{abstract}

\maketitle

\section{Introduction}
In this paper we investigate the reduced Kronecker product $[\l]_\bullet\star[\m]_\bullet$ which can be viewed roughly as the Kronecker product $[(n-\abs{\l},\l)][(n-\abs{\m},\m)]$ for $n$ big enough. Here and in the following the Kronecker product is the Kronecker product of the characters of the symmetric group over the field $\mathbb{C}$. 

First results about the reduced Kronecker product can  already be found in Murnaghan's work \cite{Mur0} of the year 1937. The reduced Kronecker product has, as we will see, some nicer properties than the ordinary Kronecker product and there is also a relation to the ordinary outer product of two irreducible characters, involving the famous Littlewood Richardson coefficients.

In the following we classify the reduced Kronecker products which are multiplicity free and those which contain less than $10$ components. Furthermore, we give general lower bounds for the number of constituents and components of a given reduced Kronecker product. We also give a lower bound for the number of pairs of components whose corresponding partitions differ by one box.

We also argue that equality of two reduced Kronecker products is only possible in the trivial case that the factors of the product are the same.

\section{Notation and Preliminaries}
We mostly follow the standard notation in \cite{Sag} or \cite{Stanley}. A partition $\l=(\l_1,\l_2,\ldots,\l_l)$ is a weakly decreasing sequence of non-negative integers where only finitely many of the $\l_i$ are positive. We regard two partitions as the same if they differ only by the number of trailing zeros and call the positive $\l_i$ the parts of $\l$. The length is the number of positive parts and we write $l(\l)=l$ for the length and $\abs{\l}=\sum_i \l_i$ for the sum of the parts. For a partition $\l$ we set $\l[n]=(n-\abs{\l},\l)$ which is again a partition for $n\geq \abs{\l}+\l_1$. Furthermore, we denote by $dp(\l)$ the number of different parts of $\l$.

We will use $\d_n$ always to refer to the staircase partition: 
\[\d_n=(n,n-1,n-2,\ldots,2,1).\]

With a partition $\l$ we associate a diagram, which we also denote by $\l$, containing $\l_i$ left-justified boxes in the $i$-th row and we use matrix-style coordinates to refer to the boxes.

The conjugate $\l^c$ of $\l$ is the diagram which has $\l_i$ boxes in the $i$-th column.

The partition $\l+(1^a)=(\l_1+1,\l_2+1,\ldots \l_a+1,\l_{a+1},\ldots)$ is obtained from $\l$ by inserting a column containing $a$ boxes.  The partition $\l\cup (a)$ is obtained from $\l$ by inserting a row containing $a$ boxes. Inserting a row is the conjugate of inserting a column: $(\l+(1^a) )^c=\l^c\cup(a)$. Both operation can be generalized to arbitrary partitions: $\l+\l'$ and $\l\cup\l'$.

For example we have:
\begin{align*}
 \yng(4,4,2,2)+\young(XXX,XX,X)&=\young(\h\h\X\h\h\X\X,\h\h\X\h\h\X,\h\h\X,\h\h), & \yng(4,4,2,2) \cup \young(XXX,XX,X)&=\young(\h\h\h\h,\h\h\h\h,XXX,\h\h,\h\h,XX,X).
\end{align*}

\begin{Def}
 We say that a partition $\l$ is larger than $\l'$ if $\l$ can be obtained from $\l'$ by repeatedly using the operations $+,\cup$ with arbitrary partitions in any order.
\end{Def}

The irreducible characters $[\l]$ of the symmetric group $S_n$ are naturally labeled by partitions $\l\vdash n$.

Since the characters are functions we can define the ordinary product $[\l][\m](\pi)=[\l](\pi)[\m](\pi)$ which is again a character and called Kronecker product. Since it is a character we can write it as a sum of irreducible characters
\[ [\l][\m]=\sum_\n g(\l,\m,\n) [\n] \] where the $g(\l,\m,\n)$ are the Kronecker coefficients. The Kronecker coefficient stays constant if one permutes the $\l,\m,\n$ and from multiplication with the sign character it follows:
\[g(\l,\m,\n)=g(\l^c,\m,\n^c).\]

\begin{Def}
We say that a character $\chi=\sum_\l c_\l [\l]$ is of cc-type $(a,b)$ if $\chi$ has $a=\sum_{c_\l\neq0} 1$ components and $b=\sum_\l c_\l$ constituents.
\end{Def}

The Kronecker product has a nice stabilizing property which we will demonstrate with an example:
\begin{align*}
[21^2][22]&=[31]\phantom{+[62]}\;\:+\phantom{2}[21^2],\phantom{+[53]+2[521]+[51^3]+[431]+[421^2]}\\
[31^2][32]&=[41]+[32]+2[31^2]\phantom{+[53]}\;\:+\phantom{2}[221]+[21^3],\phantom{+[431]+[421^2]}\\
[41^2][42]&=[51]+[42]+2[41^2]+[33]+2[321]+[31^3]\phantom{+[431]}\;\:+[221^2],\\
[51^2][52]&=[61]+[52]+2[51^2]+[43]+2[421]+[41^3]+[331]+[321^2],\\
[61^2][62]&=[71]+[62]+2[61^2]+[53]+2[521]+[51^3]+[431]+[421^2],\\
[\ast1^2][\ast2]&=[\ast1]+[\ast2]+2[\ast1^2]+[\ast3]+2[\ast21]+[\ast1^3]+[\ast31]+[\ast21^2].
\end{align*}

Murnaghan showed that this product always stabilizes if one increases the first part. For a partition $\l$ we define the sequence $\l[n]=(n-\abs{\l},\l_1,\l_2,\ldots)$ which is a weak composition of $n$ for $n\geq\abs{\l}$ and a partition of $n$ if $n\geq\abs{\l}+\l_1$.

The stabilizing property of the Kronecker product then means that for all partitions $\l,\m$ there is an integer $n$ such that for all $k\geq0$ and all partitions $\n$ we have:
 \begin{align}\label{prod}
 g(\l[n],\m[n],\n[n]) =g(\l[n+k],\m[n+k],\n[n+k]).
\end{align}

So it makes sense to define the reduced Kronecker coefficients \[\bg(\l,\m,\n)=\lim_{n \rightarrow \infty} g(\l[n],\m[n],\n[n])\] which is by Murnaghan's Theorem well defined.

There is also a correspondence of the Schur functions $s_\l$ to the irreducible character $[\l]$. It is well known (see for example \cite{Sag}) that the Schur functions can be expressed via the Jacobi-Trudi determinant
\[s_\l=\det(h_{\l_i}+i-j)_{1\leq i,j \leq n} \]
which allows to define the Schur function $s_{\l[n]}$ even for the case that $\l[n]$ is not a partition. Nevertheless we have either $s_{\l[n]}=\pm s_{\l'}$ for some partition $\l'$ or $s_{\l[n]}=0$. We set $[\l[n]]=\pm [\l']$ if $s_{\l[n]}=\pm s_{\l'}$ and $[\l[n]]=0$ if $s_{\l[n]}=0$.

For the ordinary product of Schur functions we have $s_\m s_\n=\sum_\l c(\l;\m,\n) s_\l$ with $c(\l;\m,\n)$ the famous Littlewood Richardson coefficients. The corresponding product of irreducible characters $[\m]\otimes[\n]$ is obtained by induction and referred to as outer product:
\[ [\m]\otimes[\n]:= ([\m]\times[\n])\uparrow_{S_m\times S_n}^{S_{m+n}} =\sum_\l c(\l;\m,\n) [\l].\]

Murnaghan's Theorem is then the following (see also~\cite{Thi}):
\begin{Sa}[Murnaghan,\cite{Mur0,Mur1,Mur2}]\label{Th:Mur}
For all $n\geq 0$ we have:

\[ [\l[n]] [\m[n]]=\sum_\n \bg(\l,\m,\n) [\n[n]]. \]
\end{Sa}

At first, this result seems to be surprising because we see that the decomposition of $[21^2][22]$ and $[\ast1^2][\ast2]$ are different. But we have by the Jacobi-Trudi determinant:
\begin{align*}
  [13]& =-[22], &[121]&=0, & [031]&=-[211], & [0211]&=-[1111],
\end{align*}
and so:
\begin{align*}
 [21^2][22]&=[31]+[22]+2[21^2]+[13]+2[121]+[1^4]+[031]+[021^2] \\
	   &=[31]+[22]+2[21^2]-[22]+0+[1^4]-[21^2]-[1^4]\\
	   &=[31]+[21^2].
\end{align*}

Murnaghan's Theorem inspires the following. Let $\{[\l]_\bullet| \l \textnormal{ a partition}\}$ be a basis for a $\mathbb{C}$ vector space. We then define a product on these basis vectors which can be linearly extended to the full vector space:
\[ [\l]_\bullet \star [\m]_\bullet:= \sum_\n \bg(\l,\m,\n) [\n]_\bullet. \]
We call this the reduced Kronecker product. The connection to the usual Kronecker product should be obvious by Theorem~\ref{Th:Mur}. In fact one could set $[\l]_\bullet$ formally as $[\l]_\bullet=\sum_{n\in\mathbb{N}} [\l[n]]$ with $\star$ the usual Kronecker product (this works because $[\l[n]][\m[m]]=0$ for $n\neq m$). We will call these $[\l]_\bullet$ irreducible characters even if they are not traces of representations.

In this notation our example from above simply reads:
\begin{align*}
[1^2]_\bullet \star [2]_\bullet=&[1]_\bullet +[2]_\bullet +2[1^2]_\bullet +[3]_\bullet +2[21]_\bullet +[1^3]_\bullet +[31]_\bullet +[21^2]_\bullet.
\end{align*}
This particular product can already be found in~\cite[Number 20]{Mur1}.

The reduced Kronecker coefficients has some nice properties. If $\abs{\l}=\abs{\m}+\abs{\n}$ then $\bg(\l,\m,\n)=c(\l;\m,\n)$ with $c(\l;\m,\n)$ still the LR coefficients. In \cite{BOR} Briand et al.\ showed the following:
\begin{Sa}[{\cite[Theorem 1.2]{BOR}}] \label{Sa:redBOR}
For arbitrary partitions $\l,\m$ $n=\abs{\l}+\abs{\m}+\l_1+\m_1$ is the smallest $n$ which can be chosen for equation~\eqref{prod}.

This means that the coefficients in $[\l[n]][\m[n]]$ stabilize for $n=\abs{\l}+\abs{\m}+\l_1+\m_1$.
\end{Sa}

This means that we can calculate $[\l]_\bullet \star[\m]_\bullet$ with Stembridge's Maple Package \cite{stemmaple} by calculating $[\l[n]][\m[n]]$  for $n=\abs{\l}+\abs{\m}+\l_1+\m_1$.

Furthermore, it is conjectured by Klyachko~\cite[Conjecture 6.2.4]{Klyachko} and Kirillov~\cite[Conjecture 2.33]{Kirillov} that the reduced Kronecker coefficients satisfy, like the LR coefficients, the Saturation property. So it is conjectured that if $\bg(n\l,n\m,n\m)\neq0$ for some $n\geq 1$ then $\bg(\l,\m,\m)\neq0$. Obviously this property holds if $\abs{\l}=\abs{\m}+\abs{\n}$. There are many examples which show that the ordinary Kronecker coefficients do not satisfy the saturation conjecture.

We will use the following lemma in our later proofs. It can be found in \cite[page 1098]{Mur3} and \cite[page217]{Thi}:

\begin{Le}\label{Le:n-11}
 \[ [1]_\bullet \star [\l]_\bullet= dp(\l)[\l]_\bullet +\sum_\m [\m]_\bullet \]
where the sum is over all partitions $\m$ different from $\l$ which can be obtained from $\l$ by adding a box, deleting a box or first  deleting and then adding a box.

In particular $[\l]_\bullet$ has multiplicity $dp(\l)$ and all other characters have multipli\-ci\-ty $0$ or $1$.
\end{Le}

The formulas for the ordinary Kronecker product is also easy to prove and can be found for example in \cite[Lemma~4.1]{BK}:

\begin{Le}
Let $n\geq2$ and $\l\vdash n$. Then:

\[ [n-1,1][\l] = (dp(\l)-1)[\l] + \sum_\m[\m] \]

where the sum is over all partitions $\m$ different from $\l$ which can be obtained from $\l$ by first removing and then adding a box.
\end{Le}

In \cite{Manivel} Manivel states the following lemma in the proof of his Theorem~1. 

\begin{Le}[{\cite{Manivel}}] \label{Le:addkron}
 Let both $g(\l,\m,\n),g(\l',\m',\n')\neq 0$.

Then \[g(\l+\l',\m+\m',\n+\n') \geq g(\l,\m,\n)\]
and by conjugation: \[g(\l\cup\l',\m+\m',\n\cup\n') \geq g(\l,\m,\n).\]
\end{Le}

\begin{Bem}
 Christandl et al. proved in \cite[Theorem~3.1]{CHM} that $g(\l+\l',\m+\m',\n+\n')\neq0$ for $g(\l,\m,\n),g(\l',\m',\n')\neq 0$.
\end{Bem}

It is easy to show using this results that the reduced Kronecker coefficients also obey the same property:

\begin{Le} \label{Le:addredkron}
Let $m\in\mathbb{N}$ such that $g(\l'[m],\m'[m],\n'[m])\neq0$.

Then \[\bg(\l+\l',\m+\m',\n+\n') \geq \bg(\l,\m,\n)\]
and \[\bg(\l\cup(\l'[m]),\m+\m',\n\cup(\n'[m])) \geq \bg(\l,\m,\n).\]
\end{Le}
\begin{proof}
We have
\begin{align*}
\bg(\l+\l',\m+\m',\n+\n')&=g((\l+\l')[n+m],(\m+\m')[n+m],(\n+\n')[n+m])\\
			 &=g(\l[n]+\l'[m],\m[n]+\m'[m],\n[n]+\n'[m])\\
			 &\geq g(\l[n],\m[n],\n[n])=\bg(\l,\m,\n)
\end{align*}
for $n$ large enough and by Lemma~\ref{Le:addkron}.

For $n$ large enough we also get by using Lemma~\ref{Le:addkron}:
\begin{align*}
\bg(\l\cup(\l'[m]),\m+\m',\n\cup(\n'[m]))&=g(\l[n]\cup\l'[m],\m[n]+\m'[m],\n[n]\cup\n'[m])
		 \\& \geq g(\l[n],\m[n],\n[n])= \bg(\l,\m,\n).
\end{align*}
\end{proof}

\begin{Bem}\label{Bem:main}
Let $\l'\vdash n$, then obviously we have $\bg(\l',\emptyset,\l')=1,g(\l',(n),\l')=1$.
So we have $\bg(\l+\l',\m,\n+\l')\geq\bg(\l,\m,\n)$ and also $\bg(\l\cup\l',\m,\n\cup\l')\geq\bg(\l,\m,\n)$.

This property allows us to get informations about the product $[\l]_\bullet \star[\m]_\bullet $ if we know the product $[\l']_\bullet \star[\m']_\bullet $ and $\l$ is larger than $\l'$ and  $\m$ is larger than $\m'$.

This is not the case for the ordinary Kronecker product! For example let $\l'=\m'=(3,2,1)=\yng(3,2,1)$, $\l=(4,2,1,1)=\yng(4,2,1,1)$ and $\m=(4,3,1)=\yng(4,3,1)$.

So we have $\l=\l'+(1)\cup(1)$ (in this particular case $+$ and $\cup$ commute which is not true in general) and $\m=\m'+(1^2)$.
So we know by applying Lemma~\ref{Le:addredkron} three times that if $[\n]_\bullet $ appears in $[\l']_\bullet \star[\m']_\bullet $ with multiplicity $c$ then $[(\n+(1^2)+(1) )\cup(1)]_\bullet $ appears in $[\l]_\bullet \star[\m]_\bullet $ with multiplicity at least $c$:
\begin{align*}
\bg(\l',\m',\n)&\leq \bg(\l'+(1),\m',\n+(1))\\
&\leq \bg(\l'+(1),\underbrace{\m'+(1^2)}_{\m},\n+\underbrace{(1)+(1^2)}_{(2,1)}) \\
&\leq \bg(\underbrace{\l'+(1)\cup(1)}_{\l},\m,(\n+(2,1))\cup(1)).
\end{align*}

On the other hand, the only way to get $\m$ from $\m'$ is by adding $+(1^2)$ and this has to be done directly, not in two steps. So in our case with $\m=(4,3,1)=\m'+(1^2)$ we can deduce something about $[\l][\m]$ from $[\l'][\m']$ using Lemma~\ref{Le:addkron} only if $\l\in\{\l'+(2),\l'\cup(2),\l'+(1^2),\l'\cup(1^2)\}$ which is not the case.
\end{Bem}

\begin{Bem}\label{Bem:addcc-typeRK}
We have $\l+\n\neq\m+\n$ and $\l\cup\n\neq\m\cup\n$ for $\l\neq\m$ and arbitrary $\n$. So if $[\l']_\bullet\star[\m']_\bullet$ has cc-type $(a',b')$ and $\l$ is larger than $\l'$ and $\m$ is larger than $\m'$. If $[\l]_\bullet\star[\m]_\bullet$ has cc-type $(a,b)$ then it is $a\geq a'$ and $b \geq b'$.
\end{Bem}

\begin{Bem}
In \cite[Theorem~3.1]{Gut} we showed that the Littlewood Richardson coefficients obey the same property: For $c(\l;\m,\n),c(\l';\m',\n')\neq 0$ we have \[c(\l+\l';\m+\m',\n+\n') \geq c(\l;\m,\n).\]
\end{Bem}

We let $p_n$ denote the number of partitions of $n$, $f_n$ denote the number of standard Young tableaux of size $n$ and $g_n$ denote the number of pairs of partitions of $n$ which differ by one box. 

\begin{Bem}\label{Bem:firstgpf}
 In the On-Line Encyclopedia of Integer Sequences~\cite{fcOEIS} $g_{n+2}=\bp_{n}$ has the id: A000097, $p_n$ has the id: A000041 and $f_n$ has the id: A000085. Their first terms are:
\[
\begin{array}{l|ccc|ccc|ccc|ccc|c}
n:  &1 &2 &3 &4 &5 &6 &7 &8 &9 &10 &11 &12 &13 \\
\hline
g_n:  &0 & 1 &2 &5 &9 &17 &28 & 47&73 & 114 &170 &253 &365 \\
p_n: &1 & 2 &3 &5 & 7 &11 &15 &22  &30  &42 &56  &77  &101  \\
f_n:   &1 &2 &4 &10 &26 &76 &232 &764 &2620 &9496 &35696 &140152 &568504
\end{array}\]
\end{Bem}

We will later use the following result proved in~\cite{Gut2}.
\begin{Le}[{\cite[Lemma~3.12]{Gut2}}]\label{Le:prodccpairs}
Let $\a,\b$ be partitions with $dp(\a)\geq dp(\b)=n$. Then $[\a]\otimes[\b]$ has cc-type at least $(p_{n+1},f_{n+1})$ and contains $g_{n+1}$ pairs of characters $([\n^1],[\n^2])$ such that their corresponding partitions differ only by one box.
\end{Le}

\section{Results}\label{sec:RKmain}
Using Lemma~\ref{Le:addredkron} and Remark~\ref{Bem:addcc-typeRK} it is now easy to classify the multiplicity free reduced Kronecker products and the reduced Kronecker products which contain only few components.

\begin{Sa}\label{Sa:redKronmf}
 Let $\l,\m$ be partitions. Then $[\l]_\bullet \star[\m]_\bullet $ is multiplicity free if and only if up to exchanging $\l$ and $\m$ we are in one of the following two cases:
\begin{enumerate}
 \item $\l=\emptyset,\quad \m$ arbitrary
 \item $\l=(1),\quad \m=(\a^a)$ is a rectangle
\end{enumerate}
\end{Sa}

\begin{proof}

Obviously for $\l=\emptyset$ we have the trivial product $[\emptyset]_\bullet \star[\m]_\bullet =[\m]_\bullet $ which is multiplicity free.

By Lemma~\ref{Le:n-11} $[1]_\bullet \star[\m]_\bullet $ is multiplicity free if and only if $dp(\m)\leq 1$, which is the case only for $\m$ a rectangle.

So now suppose that neither $\l$ nor $\m$ is $(1)$. Then $\l$ and $\m$ are larger than $(2)$ or $(1^2)$. So by Lemma~\ref{Le:addredkron} and Remark~\ref{Bem:main} it is enough to check if the following products have multiplicity: $[2]_\bullet \star[2]_\bullet ,[2]_\bullet \star[1^2]_\bullet $ and $[1^2]_\bullet \star[1^2]_\bullet $. We already know that $\bg(1^2,2,1^2)=2$ from our example  and so both $[2]_\bullet \star[1^2]_\bullet $ and $[1^2]_\bullet \star[1^2]_\bullet$ have multiplicity. Using Stembridge's Maple package to calculate $[2[n]][2[n]]$ for $n=8$ gives $\bg(2,2,2)=2$ which proves that also $[2]_\bullet \star[2]_\bullet$ has multiplicity.
\end{proof}

\begin{Sa}\label{Sa:redKronfc}
 Let $\l,\m$ be partitions. Then $[\l]_\bullet \star[\m]_\bullet $ has less than $10$ components if and only if up to exchanging $\l$ and $\m$ we are in one of the following cases:
\begin{enumerate}
 \item $\l=\emptyset$, $\m$ arbitrary (cc-type $(1,1)$)
 \item $\l=(1)$, $\m=(\a^a)$ is rectangle. In this case we have
	\begin{enumerate}
	 \item cc-type $(4,4)$ if $\m=(1)$
	 \item cc-type $(5,5)$ if either $\a=1$ or $a=1$
	 \item cc-type $(6,6)$ if $\a,a\geq 2$
	\end{enumerate}
 \item $\l=(1)$, $\m=(\a^a,\b^b)$ is a fat hook and we have, furthermore,
	\begin{enumerate}
	 \item $\m=(2,1)$ (cc-type $(8,9)$)
 	 \item $3$ of the values $\a-\b,\b,a,b$ are equal to $1$ (cc-type $(9,10)$)
	\end{enumerate}
 \item $\l=(1^2)$, $\m=(2)$ (cc-type $(8,10)$)
\end{enumerate}
\end{Sa}

\begin{proof}
We have already seen that $[1^2]_\bullet\star[2]_\bullet$ has the given cc-type.

To check that the other cases are true one simply uses the known formula for $[1]_\bullet\star[\l]_\bullet$.

For $\l=(1),\m=(\a^a)$ we have by Lemma~\ref{Le:n-11}:
\begin{align*}
[1]_\bullet\star[\a^a]_\bullet=&[\a^a]_\bullet+[\a^{a-1}]_\bullet+[\a^a,1]_\bullet+[\a^{a-1},\a-1]_\bullet\\&+[\a+1,\a^{a-2},\a-1]_\bullet +[\a^{a-1},\a-1,1]_\bullet
\end{align*}
where the last component appears only for $\a\geq2$ and the penultimate only for $a\geq2$.

For $\l=(1),\m=(\a^a,\b^b)$ we have by Lemma~\ref{Le:n-11}:
\begin{align*}
 [1]_\bullet \star[\a^a,\b^b]_\bullet =& 2[\a^a,\b^b]_\bullet +[\a+1,\a^{a-1},\b^b]_\bullet +[\a^a,\b+1,\b^{b-1}]_\bullet +[\a^a,\b^b,1]_\bullet  \\
 		&	+[\a+1,\a^{a-2},\a-1,\b^b]_\bullet +[\a^{a-1},\a-1,\b+1,\b^{b-1}]_\bullet  \\
		& 	+[\a^{a-1},\a-1,\b^b,1]_\bullet 	+[\a+1,\a^{a-1},\b^{b-1},\b-1]_\bullet  \\
		&	+[\a^a,\b+1,\b^{b-2},\b-1]_\bullet +[\a^a,\b^{b-1},\b-1,1]_\bullet  \\
		&	+[\a^{a-1},\a-1,\b^b]_\bullet +[\a^a,\b^{b-1},\b-1]_\bullet
\end{align*}
where the $5$th characters appears only for $a\geq2$, the $6$th only for $\a-\b\geq2$, the $9$th only for $b\geq2$ and the $10$th only for $\b\geq2$.

We will check now that all other products have at least $10$ components. From the formula above we can see that if $\l=(1)$ and $\m=(\a^a,\b^b)$ and none of the additional conditions given in the theorem is satisfied we have at least $10$ components.

Suppose now that $\l=(1)$ and $dp(\m)\geq 3$. Then $\m$ is larger than or equal to $(3,2,1)$. Since $[1]_\bullet \star[321]_\bullet $ has cc-type $(14,16)$ we know by Lemma~\ref{Le:addredkron} and Remark~\ref{Bem:addcc-typeRK} that $[1]_\bullet \star[\m]_\bullet $ has at least $14$ components.

So suppose now that both $\l,\m$ are different from $(1)$ and we are not in the situation $[1^2]_\bullet \star[2]_\bullet $.

If both $\l_1,\m_1\geq 2$ then both $\l,\m$ are larger than or equal to $(2)$. It is
\begin{align*}
[2]_\bullet \star[2]_\bullet =&[\emptyset]_\bullet +[1]_\bullet +2[2]_\bullet +[1^2]_\bullet +[3]_\bullet +2[21]_\bullet +[1^3]_\bullet +[4]_\bullet +[31]_\bullet +[2^2]_\bullet
\end{align*} and so  $[2]_\bullet \star[2]_\bullet $ has cc-type $(10,12)$ (this product can already be found in \cite[Number 19]{Mur1}). It follows by Remark~\ref{Bem:addcc-typeRK} that $[\l]_\bullet \star[\m]_\bullet $ also has at least $10$ components.

If both $l(\l),l(\m)\geq 2$ then both $\l,\m$ are larger than or equal to $(1^2)$.
It is \begin{align*}
[1^2]_\bullet \star[1^2]_\bullet =&[\emptyset]_\bullet +[1]_\bullet +2[2]_\bullet +[1^2]_\bullet +[3]_\bullet \\&+2[21]_\bullet +[1^3]_\bullet +[2^2]_\bullet +[21^2]_\bullet +[1^4]_\bullet
\end{align*} and so $[1^2]_\bullet \star[1^2]_\bullet $ also has cc-type $(10,12)$ (this product can already be found in \cite[Number 29]{Mur1}). It follows by Remark~\ref{Bem:addcc-typeRK} that $[\l]_\bullet \star[\m]_\bullet $ also has at least $10$ components.

So we may now suppose that $\l=(\l_1)$ and $\m=(1^m)$ with $\l_1,m\geq2$ and at least one of $\l_1,m\geq3$. Using Stembridge's Maple package~\cite{stemmaple} we check that $[2]_\bullet \star[1^3]_\bullet $ has cc-type $(10,13)$ and $[3]_\bullet \star[1^2]_\bullet $ has cc-type $(11,13)$ (these products can also be found in \cite[Numbers $23$ and $30$]{Mur1}). It follows by Remark~\ref{Bem:addcc-typeRK} that $[\l]_\bullet \star[\m]_\bullet $ also has at least $10$ components.
\end{proof}

In the reduced Kronecker products appear characters whose corresponding partitions are of different size. So in the following if we say that two partitions $\l,\m$ differ by one box we mean that $\abs{\l\cap\m}=\max(\abs{\l},\abs{\m})-1$. So $\m$ can be obtained from $\l$ by deleting a box, adding a box or first deleting and then adding a box.

\begin{Sa}\label{Sa:redKronpair}
Let $\l,\m$ be partitions with $dp(\l)=n\geq dp(\m)=m\geq1$. Then $[\l]_\bullet\star[\m]_\bullet$ contains  at least
\begin{itemize}
 \item $n^2+1+\max(p_{m+1},n+1)$ components,
 \item $n^2+n +\max(f_{m+1},n+1)$ constituents,
 \item $n^3+n+1+\max\left(g_{m+1},\frac{1}{2}(n+1)n\right)$ pairs of characters $([\n^1]_\bullet,[\n^2]_\bullet)$ such that their corresponding partitions $\n^1,\n^2$ differ only by one box.
\end{itemize}
 
\end{Sa}
\begin{proof}
We first investigate the product $[1]_\bullet\star[\d_n]_\bullet$ which is well known by Lemma~\ref{Le:n-11}. We have:
 \[ [1]_\bullet \star [\d_n]_\bullet= n[\d_n]_\bullet +\sum_\n [\n]_\bullet\]
where the sum is over all partitions $\n$ different from $\d_n$ which can be obtained from $\d_n$ by adding a box, deleting a box or first  deleting and then adding a box.

We label the characters  appearing in $[1]_\bullet \star [\d_n]_\bullet$ by the four ways in which they can be obtained. Let $a$ denote $[\d_n]_\bullet$, $b$ label those characters whose partitions are obtained from $\d_n$ by adding a box, $c$ those which are obtained by deleting a box and $d$ those which are obtained by first deleting and then adding a box.

So we have $1$ character labeled $a$ with multiplicity $n$, $n+1$ labeled $b$, $n$ labeled $c$ and $n(n-1)$ labeled $d$.

Now $\l$ is larger than $\d_n$ and $\m$ is larger than $(1)$. Using Remarks~\ref{Bem:main} and \ref{Bem:addcc-typeRK} we know that $[\l]_\bullet\star[\m]_\bullet$ has at least $1+n+1+n+n(n-1)=n^2+1+(n+1)$ components and $n+n+1+n+n(n-1)=n^2+n+(n+1)$ constituents. The partitions of characters labeled $b$ have size $\abs{\d_n}+\abs{(1)}$. So in $[\l]_\bullet\star[\m]_\bullet$ there are also at least $n+1$ characters whose corresponding partitions are of size $\abs{\l}+\abs{\m}$. But we have also $\bg(\l,\m,\n)=c(\n;\m,\l)$ for $\abs{\n}=\abs{\m}+\abs{\l}$ so we can deduce from  Lemma~\ref{Le:prodccpairs} that there are also at least $p_{m+1}$ components and $f_{m+1}$ constituents in $[\l]_\bullet\star[\m]_\bullet$ whose corresponding partitions are of size $\abs{\m}+\abs{\l}$. So we have in total at least $n^2+1+\max(p_{m+1},n+1)$ components and $n^2+n +\max(f_{m+1},n+1)$ constituents.

We now prove the lower bound for the number of pairs of partitions which differ by one box.

We  have in $[1]_\bullet\star[\d_n]_\bullet$ the following number of pairs of characters whose diagrams differ by one box (where type $(a,b)$ means that one character is labeled $a$ and the other $b$):
\[
\begin{array}{r|cc|cc|c}
\textnormal{Type}:  & (a,a) & (a,b) & (a,c) & (a,d) & ( b,b)    \\
\textnormal{Number}: & 0     & n+1   & n     & n(n-1) & \frac{1}{2}(n+1)(n)\\
\hline
\textnormal{Type}:  & (b,c) & (b,d) & (c,c)                    & (c,d) & (d,d)    \\
\textnormal{Number}: &  0    & n(n-1) & \frac{1}{2}n(n-1) & n(n-1) &n(n-1)(n-\frac{5}{2}).
\end{array}\]

The first four numbers are clear since all diagrams labeled $b,c$ or $d$ differ from $\d_n$ by only one box. The numbers for type $(b,b)$ and type $(c,c)$ are also clear since all diagrams of characters of type $b$ differ from one another only by the additional box. The same goes for all characters of type $c$. Furthermore, the  diagrams of characters of type $b$ have two more boxes than those of type $c$ hence the $0$ for type $(b,c)$.

We  now look at pairs of type $(b,d)$. Suppose the partition $\a$ corresponds to a character of type $b$ and $\b$ to a character of type $d$ and $\a$ and $\b$ differ by only one box. If $\a=\d_n+(1)$ (so the additional box of $\a$ is in the first row) then there are  $n-1$ partitions $\b$ with character of type $d$ ($\b$ is  obtained from $\a$ by deleting a box in a row other than the first row). The same goes if $\a=\d_n \cup(1)$. So now suppose $\a$ is neither $\d_n+(1)$ nor $\d_n\cup(1)$. We have $n-1$ of those $\a$ and in each case there are  $n-2$ $\b$ whose corresponding character is of type $d$  such that $\a$ and $\b$ differ by one box. So we have in total:
\[2(n-1)+(n-1)(n-2)=n(n-1).\]

For type $(c,d)$ we can choose one of the $n$ characters of type $c$ say with corresponding partition $\a$. There are $n-1$ places in which we can add a box to $\a$ to obtain a partition $\b$ whose corresponding character is of type $d$. So we have $n(n-1)$ of those pairs.

Let us now check the number of pairs of type $(d,d)$. Let $\a$ and $\b$ be partitions with  corresponding character of type $d$ which differ by one box. For $\a$ we can choose any of the $n(n-1)$ partitions. Suppose $\a$ is obtained from $\d_n$ by deleting the box $A$ and then adding the box $B$. For $\b$ we can then choose any of the $n-2$ partitions which are obtained from $\d_n$ by also deleting $A$ and then adding a  box different from $B$ or any of the $n-3$ partitions which are obtained from $\d_n$ by deleting a box different from $A$ (such that the box $B$ can be added afterwards) and then adding the box $B$.
Since we count each pair only once we get a factor $\frac{1}{2}$ and so in total:
\[\frac{1}{2}n(n-1)(n-2+n-3)=n(n-1)(n-\frac{5}{2}).\]
Adding all the number of pairs except pairs of type $(b,b)$ we get:
\begin{align*}
 &n+1 + n + n^2-n+ n^2 -n +\frac{1}{2}n^2-\frac{1}{2}n+n^2-n +n^3-\frac{7}{2}n^2+\frac{5}{2}n=\\
&n^3+n+1.
\end{align*}
Now $\l$ is still larger than $\d_n$ and $\m$ is larger than $(1)$. Notice that if $\n^1$ and $\n^2$ differ by one box then for every partition $\a$ also $\n^1+\a$ and $\n^2+\a$ (resp.\ $\n^1\cup\a$ and $\n^2\cup\a$) also differ by one box.  Using this fact and again Remarks~\ref{Bem:main} and \ref{Bem:addcc-typeRK}  different pairs of characters whose corresponding partitions differ by one box in $[1]_\bullet\star[\d_n]_\bullet$ give different pairs in $[\m]_\bullet\star[\l]_\bullet$. Notice that the pairs of type $(b,b)$ correspond to pairs of characters whose partitions have as size the sum of the sizes of $(1)$ and $\d_n$. So the number of pairs of characters in $[\m]_\bullet\star[\l]_\bullet$ whose partitions are of size $\abs{\m}+\abs{\l}$ and which differ by one box is at least $\frac{1}{2}n(n-1)$, the number of pairs of type $(b,b)$. But since we have $\bg(\l,\m,\n)=c(\n;\m,\l)$ if $\abs{\n}=\abs{\m}+\abs{\l}$ we know by Lemma~\ref{Le:prodccpairs} that the number of those pairs is also at least $g_{m+1}$.

This gives in total
 \[n^3+n+1+\max\left(g_{m+1},\frac{1}{2}(n+1)n\right) \]
pairs of characters whose corresponding partitions differ by only one box.
\end{proof}

Reiner et al.\ proved in \cite[Section 6]{hlRSW} the following, using the Jacobi-Trudi determinant but no LR combinatorics:

\begin{Le}[{\cite[Section 6]{hlRSW}}]\label{RSW}
 Let $[\l]\otimes[\m]=[\l']\otimes[\m']$.

Then either $\l=\l',\m=\m'$ or $\l=\m',\m=\l'$.
\end{Le}

The used property $\bg(\l,\m,\n)=c(\n;\m,\l)$ if $\abs{\n}=\abs{\m}+\abs{\l}$ answers also the question of equality of reduced Kronecker products:

\begin{Le}\label{Le:redKroneq}
 Let $\l,\l',\m,\m'$ be partitions with $[\l]_\bullet\star[\m]_\bullet=[\l']_\bullet\star[\m']_\bullet$.

Then either $\l=\l',\m=\m'$ or $\l=\m',\m=\l'$.
\end{Le}
\begin{proof}
 This follows directly from  $\bg(\l,\m,\n)=c(\n;\m,\l)$ if $\abs{\n}=\abs{\m}+\abs{\l}$ and Lemma~\ref{RSW} which states the result of the Lemma for the outer ordinary product.
\end{proof}

\begin{Bem}
 Note again that the situation for the ordinary Kronecker product is not as nearly as nice as for the reduced Kronecker product. We can't determine anything about the number of pairs of characters in the ordinary Kronecker product whose corresponding partitions differ by only one box by analyzing the corresponding ordinary Kronecker product. We have:
\[ [\frac{n(n+1)}{2}-1,1][\d_n] = (n-1)[\d_n]+\sum_\n [\n]\]
where the sum is over all partitions $\n$ different from $\d_n$ which can be obtained from $\d_n$ by first  deleting and then adding a box. As already mentioned in Remark~\ref{Bem:main} only in some cases we get informations about the product $[\l][\m]$ even if $dp(\l)=n$ (so $\l$ is larger than $\d_n$) and $\m$ is larger than $(\frac{n(n+1)}{2}-1,1)$ from facts about the product $[\frac{n(n+1)}{2}-1,1][\d_n]$.

Furthermore, for the ordinary Kronecker product we have $[\l][\m]=[\l^c][\m^c]$ which we have already seen in the examples does not hold for the reduced Kronecker product, and by Lemma~\ref{Le:redKroneq} never holds.
\end{Bem}

{\bfseries Acknowledgement:} John Stembridge's "SF-package for maple" \cite{stemmaple} was very helpful for computing examples. Furthermore, my thanks go to Christine Bessenrodt, Emmanuel Briand, Laurent Manivel and Matthias Christandl for helpful discussions regarding Lemma~\ref{Le:addkron}. This paper was inspired by some questions of my thesis supervisor Christine Bessenrodt who in \cite{BK} conjectures also a classification of the usual Kronecker products containing only few components.

My thanks also go to Emmanuel Briand and Mercedes Rosas who organized the School and Workshop "Mathematical Foundations of Quantum Information".

\end{document}